\documentclass{amsart}

\usepackage{tikz}
\usetikzlibrary{arrows}
\usepackage{color}
\usepackage[colorlinks=true,linktoc=all,linkcolor=blue,citecolor=red,filecolor=pink,urlcolor=blue]{hyperref}
\usepackage[hyperref=true,style=numeric,backend=bibtex,maxnames=100,doi=false,isbn=false,url=false]{biblatex}
\bibliography{biblio.bib}
\usepackage{csquotes}

\usepackage{amssymb,amsmath,latexsym,amsmath}
\theoremstyle{plain}                    
\newtheorem{thm}{Theorem}[section]

\newtheorem{prop}[thm]{Proposition}
\newtheorem{cor}[thm]{Corollary}
\theoremstyle{definition}
\newtheorem{defin}[thm]{Definition}
\newtheorem{ex}[thm]{Example}
\theoremstyle{remark}
\newtheorem{rmk}[thm]{Remark}
\numberwithin{equation}{section}

\newcommand{\real}{\mathbb{R}}
\newcommand{\cmp}{\mathbb{C}}
\newcommand{\intz}{\mathbb{Z}}
\newcommand{\nat}{\mathbb{N}}

\newcommand{\hyp}{\mathbb{H}}

\newcommand{\cp}{\mathbb{C}\mathbb{P}^1}
\newcommand{\rp}{\mathbb{R}\mathbb{P}^1}
\newcommand{\pslc}{\mathrm{PSL}_2\cmp}

\newcommand{\pslr}{\mathrm{PSL}_2\real}

\newcommand{\qf}{Fuchsian }

\begin{document}
\title[MULTI(DE)GRAFTING QUASI-FUCHSIAN $\cp$-STRUCTURES VIA BUBBLES]{MULTI(DE)GRAFTING 
QUASI-FUCHSIAN COMPLEX PROJECTIVE STRUCTURES VIA BUBBLES}
\author{LORENZO RUFFONI}
\address{Dipartimento di Matematica - Universit\`a di Bologna, Piazza di Porta San Donato 5, 40126, Bologna, Italy}
\curraddr{Department of Mathematics - Florida State University, 1017 Academic Way, Tallahassee, FL 32306-4510, USA}
\email{lorenzo.ruffoni2@gmail.com}

\keywords{Complex projective structures, quasi-Fuchsian holonomy, grafting, bubbling.}
\subjclass[2010]{57M50, 20H10, 14H15}
\date{\today}

\begin{abstract}
We show that the simultaneous (de)grafting of a complex projective structure with 
quasi-Fuchsian holonomy along a multicurve  can be performed by a simple sequence of one bubbling 
and one debubbling. As a consequence we obtain that any complex projective structure with 
quasi-Fuchsian holonomy $\rho:\pi_1(S)\to$ PSL$_2\mathbb{C}$  can be joined to the corresponding 
uniformizing structure $\sigma_\rho$ by a simple sequence of one bubbling and one 
debubbling, with a stopover in the space of branched complex projective structures.
\end{abstract}

\maketitle
\tableofcontents

\section{Introduction}
Complex projective structures with (quasi-)Fuchsian holonomy arise classically in the theory of 
(simultaneous) uniformization of Riemann surfaces by means of hyperbolic metrics (see 
\cite{P},\cite{BE}). They are geometric structures locally modelled on the geometry defined on the 
Riemann sphere $\cp$ by the natural action of $\pslc$ by M\"obius transformations, and admit a 
rich deformation theory (see \cite{DU} for a survey). \par
A classical result of Goldman (see \cite{GO}) states that any complex projective 
structure with Fuchsian holonomy is obtained from a hyperbolic surface by $2\pi$-grafting along a multicurve, i.e. 
by replacing some simple closed geodesics by annuli endowed with suitable projective structures, and that a similar statement holds for quasi-Fuchsian representations. 
Motivated by the study of ODEs on Riemann surfaces, Gallo-Kapovich-Marden asked in \cite{GKM} 
whether it is possible to obtain a similar statement for the more general class of branched complex 
projective structures (in which some cone points of angle $2\pi k$ for $k\in \nat$ are allowed), 
i.e. to describe all the (branched) complex projective structures with a fixed 
holonomy representation by means of elementary geometric surgeries. The surgery they propose to 
produce ramification is called bubbling, and consists in replacing a simple arc on a surface with 
a disk with a suitable projective structure. Building on Goldman's Theorem and results by 
Calsamiglia-Deroin-Francaviglia from \cite{CDF}, it was shown in \cite{R} that in Fuchsian holonomy 
essentially every structure with two simple branch points is obtained via $2\pi$-graftings and bubblings on a 
hyperbolic surface; once again a similar statement holds for quasi-Fuchsian representations.\par 
The purpose of this paper is to investigate the relationship between these two surgeries, for 
structures with quasi-Fuchsian holonomy. The following is the main result we prove (see Theorem
\ref{severalgraftingregions} below).
\begin{thm}
Let $\rho:\pi_1(S)\to \pslc$ be quasi-Fuchsian.
 Let $\sigma_0$ be a complex projective structure with holonomy $\rho$ and $\beta \subset \sigma_0$ 
a bubbleable arc which transversely crosses any grafting annulus it meets. Then there exist a 
complex projective structure $\sigma_0'$ with the same holonomy $\rho$ and a bubbleable arc 
$\beta'\subset \sigma_0'$ which avoids all the grafting annuli of $\sigma_0'$ and such that 
$Bub(\sigma_0,\beta)=Bub(\sigma_0',\beta')$.
\end{thm}
By choosing a suitable arc one then gets that any multigrafting and any multidegrafting can be 
obtained via a simple sequence of just one bubbling and one debubbling (see 
\ref{multidegraftviabubbltohyp} below). This generalises  \cite[Theorem 5.1]{CDF}, according to 
which any simple grafting can be realised by a sequence of one bubbling and one debubbling.\par
As a consequence we deduce an explicit uniform bound on the number of surgeries generically needed 
to join a couple of branched complex projective structures with the same quasi-Fuchsian holonomy and at most 
two simple branch points (see \ref{6steps} below for a more precise statement).
\begin{cor}
 Let $\rho:\pi_1(S)\to \pslc$ be quasi-Fuchsian and let  $\sigma,\tau$ be a generic couple of 
structures with holonomy $\rho$ and at most two simple branch points. Then $\tau$ is obtained from 
$\sigma$ via a sequence of at most three bubblings and three debubblings.
\end{cor}
The above result also allows us to produce examples of branched structures which are simultaneously 
obtainable as bubbling in different ways; this phenomenon shows that branched structures do not 
have a well-defined underlying unbranched structure even in the case they are realised via bubbling 
(see \ref{nowelldefunbranched} below for more details).\par
The paper is organised as follows. First of all we recall the definitions of complex 
projective structures, of \qf representations and of the two surgeries we are interested in 
(grafting in \ref{s_grafting} and bubbling in \ref{s_bubbling}). The focus is always on unbranched 
structures, the branched ones being introduced and used only as a tool for the study of the former. 
The second part contains the proof of the above results; the strategy is  the following: we 
first perform a bubbling along an arc which crosses the grafting regions and then look for a new 
bubble in the branched structure with some specified behaviour. This procedure is first exemplified 
in the case of a crossing of a simple grafting in \ref{s_annulus}, then the case of parallel 
graftings is considered in \ref{s_region} and finally  the general case is addressed in 
\ref{s_general}. For the sake of simplicity, we will work only with Fuchsian representation, but 
everything extends automatically to the quasi-Fuchsian case via a conjugation by a quasi-conformal map.\par
\vspace{.2cm}
\noindent \textbf{Acknowledgements}. I want to thank Stefano Francaviglia for many useful and colourful discussions about bubbles.
 
\section{Complex projective structures with quasi-Fuchsian holonomy}\label{s_projstr}
Let $S$ be a closed, connected and oriented surface of genus $g\geq 2$. We denote the Riemann sphere 
by $\cp=\cmp \cup \{\infty\}$  and its group of biholomorphisms by $\pslc$. A complex projective 
structure is a ($\pslc,\cp$)-structure (see \cite{DU} for a great survey); more precisely we adopt 
the following definition.
\begin{defin}\label{def_cpstructure}
A complex projective chart on $S$ is a couple $(U,\varphi)$ where $U\subset S$ is an open 
subset of $S$ and $\varphi : U \to \varphi(U) \subseteq \cp$ is an orientation 
preserving diffeomorphism with an open subset of $\cp$.  Two charts $(U,\varphi)$ and 
$(V,\psi)$ are compatible if $\exists \ g \in \pslc$ such that $\psi=g\varphi$ on $U\cap V$. A 
\textbf{complex projective structure} $\sigma$ on $S$ is the datum 
of a maximal atlas of complex projective charts.
\end{defin}
Performing analytic continuation of local charts and local change of coordinates along paths in 
$S$, we can associate to a given structure an equivalence class of development-holonomy pairs, i.e. 
of couples ($dev,\rho$) where $dev:\widetilde{S}\to \cp$ is an orientation preserving local 
diffeomorphism (called the developing map) which is equivariant with respect to a representation 
$\rho:\pi_1(S)\to \pslc$ (called the holonomy representation); such a pair is well-defined only up 
to the $\pslc$-action $g.(dev,\rho)=(g dev,g\rho g^{-1})$.\par
Since the geometry $(\pslr,\hyp^2)$ embeds in the geometry $(\pslc,\cp)$, every hyperbolic metric 
on $S$ provides an example of complex projective structure, namely the one obtained by 
$\hyp^2/\rho(\pi_1(S))$. However a general complex projective structure is not  uniformizable, in 
the sense that the developing map fails to be a diffeomorphism onto an open domain of $\cp$. 
Hyperbolic structures (which are uniformizable as projective structures) play a special role among 
complex projective structure with the same holonomy. In this paper we are concerned with the study 
of structures whose holonomy admits such a hyperbolic structure; we adopt the following definitions.
\begin{defin}
A group $\Gamma \subset \pslc$ is said to be Fuchsian if it is conjugated to a discrete cocompact 
subgroup of $\pslr$. A group $\Gamma \subset \pslc$ is said to be quasi-Fuchsian if its action on 
$\cp$ is conjugated to that of a Fuchsian group via an orientation preserving 
homeomorphism of $\cp$.
\end{defin}
By a classical result of Bers (see \cite{BE}), if $\Gamma$ is finitely generated, then this 
homeomorphism can be chosen to be quasi-conformal.
Such a group preserves a decomposition  $\cp=\Omega^+_\Gamma \cup \Lambda_\Gamma \cup 
\Omega^-_\Gamma$ into a couple of disks $\Omega^\pm_\Gamma$ (the discontinuity domain) and a Jordan 
curve $\Lambda_\Gamma$ (the limit set); for Fuchsian groups this is just the decomposition 
$\cp=\mathcal{H}^+ \cup \Lambda_\rho \cup \mathcal{H}^-$, where $\mathcal{H}^\pm$ is the 
upper/lower-half plane in $\cmp$.
\begin{defin}
A representation $\rho:\pi_1(S)\to \pslc$ is a \textbf{quasi-Fuchsian representation} if its image 
is a quasi-Fuchsian group and if there exists an orientation preserving $\rho$-equivariant 
diffeomorphism $f:\widetilde{S}\to \Omega^+_{\rho(\pi_1(S))}$. The structure 
$\sigma_\rho=\Omega^+_{\rho(\pi_1(S))}/\rho(\pi_1(S))$ is called the uniformizing structure for 
$\rho$. When $\rho$ is Fuchsian, this is a hyperbolic structure.
\end{defin}
A diffeomorphism as in this definition is precisely a developing map for the structure 
$\sigma_\rho$. More generally we can equivariantly pullback the $\rho$-invariant decomposition of 
$\cp$ via the developing map of any quasi-Fuchsian structure $\sigma$ to obtain a decomposition of $S$.

\begin{defin}
 The \textbf{geometric decomposition} of a quasi-Fuchsian structure is the decomposition 
$S=\sigma^+\cup\sigma^\real\cup \sigma^-$, where $\sigma^\pm$ is the set of points developing to 
$\Omega^\pm_{\rho(\pi_1(S))}$ and $\sigma^\real$ is the set of points developing to 
$\Lambda_{\rho(\pi_1(S))}$. They are respectively called the positive/negative part and the real 
curve of $\sigma$; a connected component of $\sigma^+$ (respectively $\sigma^-$, $\sigma^\real$) 
will be called a positive (respectively negative, real) component of the real decomposition
\end{defin}

From now on we restrict our attention to the case of Fuchsian representations for simplicity, but all the methods and statements readily extend to the quasi-Fuchsian case just via conjugation by a quasi-conformal map. This make it easier to state that the pieces of the geometric decomposition are naturally equipped with nice geometric structures.
It follows from the definitions that $\sigma^\real$ is a finite union of simple closed curves on 
$S$ equipped with a ($\pslr,\rp$)-structure and that $\sigma^\pm$ is a finite union of subsurfaces 
endowed with hyperbolic structures;  these hyperbolic metrics are indeed complete, hence these 
pieces decompose as a union of a compact convex core and several annular ends. Each end is a 
semi-infinite embedded hyperbolic annulus whose boundary consists of one real curve 
and one geodesic (see \cite{GO} for more details). For example, the geometric decomposition for the 
uniformizing structure $\sigma_\rho$ consists of a single positive component which coincides with 
the whole surface. A couple of geometric surgeries are known in the literature, which allow to 
produce structures with more complicated geometric decompositions without changing the holonomy. We 
now introduce them, and investigate a relationship between them in the second part of the paper.

\subsection{Grafting}\label{s_grafting}

The first surgery was introduced by Maskit in \cite{MAS2} to produce examples of exotic projective 
structures, i.e. structures with surjective non injective developing map. Let $\sigma$ be a 
structure defined by a development-holonomy pair $(dev,\rho)$.
\begin{defin}
 A simple closed curve $\gamma \subset S$ is graftable with respect to $\sigma$ if $\rho(\gamma)$ is 
a non-elliptic loxodromic and $\gamma$ is injectively developed, i.e. the restriction of $dev$ to 
any of its lifts $\widetilde{\gamma}\subset \widetilde{S}$ is injective.
\end{defin}
Since $dev$ is $\rho$-equivariant, the developed image of a graftable curve is an 
embedded arc in $\cp$ joining the two fixed points of $\rho(\gamma)$; moreover $\rho(\gamma)$ acts 
freely and properly discontinuously on $\cp \setminus \overline{dev(\widetilde{\gamma})}$ and the 
quotient is an annulus endowed with a natural complex projective structure. 
\begin{defin}
 Let $\gamma\subset S$ be a graftable curve with respect to $\sigma$. For any lift 
$\widetilde{\gamma}$ of $\gamma$ we cut $\widetilde{S}$ along it and a copy of $\cp$ along 
$\overline{dev(\widetilde{\gamma})}$, and glue them together equivariantly via the developing map. 
We obtain a simply connected surface $\widetilde{S}'$ to which the action 
$\pi_1(S)\curvearrowright \widetilde{S}$ and the map $dev:\widetilde{S}\to \cp$ naturally extend, 
so that the quotient is naturally endowed with a new projective structure. We 
call this structure the \textbf{grafting} of $\sigma$ along $\gamma$ and denote it by 
$Gr(\sigma,\gamma)$. The surface $\sigma \setminus \gamma$  embeds in $Gr(\sigma, 
\gamma)$ and the complement is the annulus $A_\gamma=(\cp \setminus 
\overline{dev(\widetilde{\gamma})})/\rho(\gamma)$, which we call the grafting annulus associated to 
$\gamma$. 
\end{defin}
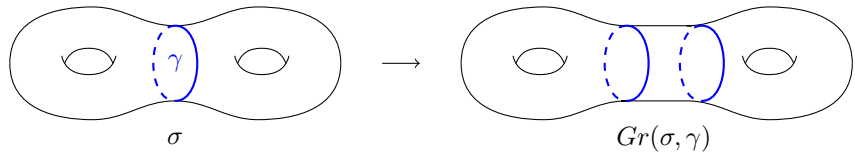
\begin{figure}[h]

\begin{center}

\begin{tikzpicture}[scale=0.5]
\node at (0,-2) {$\sigma$};
\node at (0,0) {\color{blue} $\gamma$};
\draw[xscale=0.8] (-5.5,0) to[out=90,in=180] (-2.75,1.5) to[out=0,in=180] (0,1) to[out=0,in=180] (2.75,1.5) to[out=0,in=90] (5.5,0) ;
\draw[rotate=180,xscale=0.8] (-5.5,0) to[out=90,in=180] (-2.75,1.5) to[out=0,in=180] (0,1) to[out=0,in=180] (2.75,1.5) to[out=0,in=90] (5.5,0) ;

\draw[xshift=-0.3cm,yshift=-0.1cm,xscale=0.8]  (-3.3,0.1) to[out=65,in=180] (-2.5,0.5) to[out=0,in=115] (-1.7,0.1);
\draw[xshift=-0.3cm,yshift=-0.1cm,xscale=0.8]  (-3.4,0.3) to[out=-75,in=180] (-2.5,-0.2) to[out=0,in=-105] (-1.6,0.3);
\draw[xshift=4.3cm,yshift=-0.1cm,xscale=0.8]  (-3.3,0.1) to[out=65,in=180] (-2.5,0.5) to[out=0,in=115] (-1.7,0.1);
\draw[xshift=4.3cm,yshift=-0.1cm,xscale=0.8]  (-3.4,0.3) to[out=-75,in=180] (-2.5,-0.2) to[out=0,in=-105] (-1.6,0.3);

\draw[blue,dashed,thick] (0,1) to[out=180,in=180] (0,-1);
\draw[blue,thick] (0,1) to[out=0,in=0] (0,-1);


\draw[->] (5.5,0) -- (6.5,0);

\begin{scope}[xshift=12cm]
 \node at (1,-2) {$Gr(\sigma,\gamma)$};
\draw[xscale=0.8] (-5.5,0) to[out=90,in=180] (-2.75,1.5) to[out=0,in=180] (0,1) to[out=0,in=0] (2,1) 
to[out=0,in=180] (4.75,1.5) to[out=0,in=90] (7.5,0) ;
\draw[rotate=180,xscale=0.8,xshift=-2cm] (-5.5,0) to[out=90,in=180] (-2.75,1.5) to[out=0,in=180] 
(0,1) to[out=0,in=0] (2,1) to[out=0,in=180] (4.75,1.5) to[out=0,in=90] (7.5,0) ;

\draw[xshift=-0.3cm,yshift=-0.1cm,xscale=0.8]  (-3.3,0.1) to[out=65,in=180] (-2.5,0.5) 
to[out=0,in=115] (-1.7,0.1);
\draw[xshift=-0.3cm,yshift=-0.1cm,xscale=0.8]  (-3.4,0.3) to[out=-75,in=180] (-2.5,-0.2) 
to[out=0,in=-105] (-1.6,0.3);
\draw[xshift=5.8cm,yshift=-0.1cm,xscale=0.8]  (-3.3,0.1) to[out=65,in=180] (-2.5,0.5) 
to[out=0,in=115] (-1.7,0.1);
\draw[xshift=5.8cm,yshift=-0.1cm,xscale=0.8]  (-3.4,0.3) to[out=-75,in=180] (-2.5,-0.2) 
to[out=0,in=-105] (-1.6,0.3);

\draw[blue,dashed,thick] (0,1) to[out=180,in=180] (0,-1);
\draw[blue,thick] (0,1) to[out=0,in=0] (0,-1);
\draw[blue,dashed,thick] (2,1) to[out=180,in=180] (2,-1);
\draw[blue,thick] (2,1) to[out=0,in=0] (2,-1);
\end{scope}

\end{tikzpicture}

\end{center}
\caption{Grafting a surface}
\end{figure}
This construction can of course be extended to perform simultaneous graftings on a
disjoint collection of graftable curves. It is also possible to attach an integer weight $M\in 
\nat$ to a graftable curve and to perform an $M$-fold grafting along it by gluing not just one 
copy of $\cp \setminus \overline{dev(\widetilde{\gamma})}$ but $M$ copies of it, attached in a 
chain of  length $M$ along their boundaries. The corresponding region in the surface is a chain 
 $A_\gamma=\cup_{k=1}^M A_{\gamma}^k$ of $M$ copies of the annulus $(\cp \setminus 
 \overline{dev(\widetilde{\gamma})})/\rho(\gamma)$, which we call the \textbf{grafting region} 
 associated to $M\gamma$, and we reserve the term grafting annulus for each individual 
 $A_\gamma^k$. This generalisation allows to perform a grafting along any graftable multicurve; 
 we call this operation \textbf{multigrafting}. The inverse operation is called a 
(\textbf{multi})\textbf{degrafting}. Notice that both operations preserve the holonomy of the 
structure.
 \begin{ex}\label{ex_graftgeod}
 The easiest example consists in grafting a simple geodesic on a 
 hyperbolic surface; for such a structure every simple essential curve $\gamma$ is graftable, since 
the holonomy is purely hyperbolic and the developing map is globally injective.
 \begin{figure}[h]
\begin{center}

\begin{tikzpicture}[xscale=0.75]
\draw[blue,dashed,thick] (0,1) to[out=-120,in=120] (0,-1);
\draw[blue,thick] (0,1) to[out=-60,in=60] (0,-1);
 \node at (0,-2) {\color{blue} $\gamma^-$};

 \draw[blue,dashed,thick] (6,1) to[out=-120,in=120] (6,-1);
\draw[blue,thick] (6,1) to[out=-60,in=60] (6,-1);
 \node at (6,-2) {\color{blue} $\gamma^+_R$};
 
 \draw[blue,thick] (-6,1) to[out=-120,in=120] (-6,-1);
\draw[blue,thick] (-6,1) to[out=-60,in=60] (-6,-1);
 \node at (-6,-2) {\color{blue} $\gamma^+_L$};
 
 \draw[red,dashed,thick] (3,2) to[out=-120,in=120] (3,-2);
\draw[red,thick] (3,2) to[out=-60,in=60] (3,-2);
 \node at (3,-3) {\color{red} $l_R$};

 \draw[red,dashed,thick] (-3,2) to[out=-120,in=120] (-3,-2);
\draw[red,thick] (-3,2) to[out=-60,in=60] (-3,-2);
 \node at (-3,-3) {\color{red} $l_L$};

 \draw (0,1) to[out=0,in=215] (3,2) to[out=-35,in=180] (6,1);
 \draw (0,-1) to[out=0,in=145] (3,-2) to[out=35,in=180] (6,-1);
\draw (-6,1) to[out=0,in=215] (-3,2) to[out=-35,in=180] (0,1);
 \draw (-6,-1) to[out=0,in=145] (-3,-2) to[out=35,in=180] (0,-1);
 
  \node at (-4.5,0) {$+$};
  \node at (-1.5,0) {$-$};
  \node at (1.5,0) {$-$};
  \node at (4.5,0) {$+$};

\end{tikzpicture}

\end{center}
 \caption{Geometric decomposition of the grafting annulus of 
$Gr(\sigma_\rho,\gamma)$}\label{pic_geomdecompgraftingannulus}
 \end{figure}
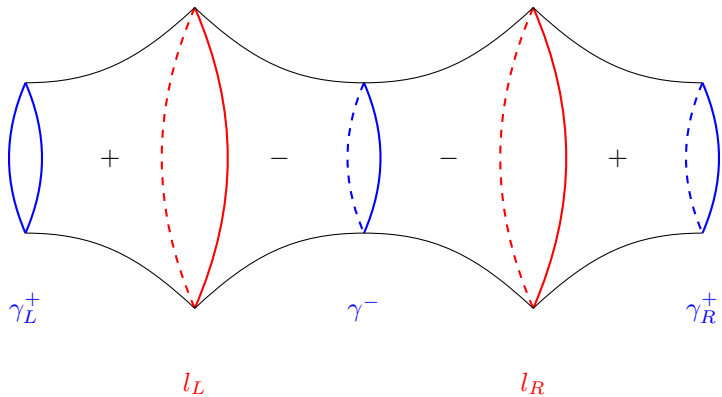
 In the geometric decomposition of the  grafting annulus $A_\gamma$ we see a 
 negative annulus (coming from the lower-half plane), bounded by a couple of real curves $l_R,l_L$ 
(coming from $\rp$) and then a couple of positive annuli (coming from the upper-half plane).
The boundary of $A_\gamma$ consists of a couple of simple closed geodesic $\gamma^+_R,\gamma^+_L$ 
coming 
from $\gamma$ and developing to the  positive part of the invariant axis of $\rho(\gamma)$, whereas 
the core of the negative annulus is a simple closed geodesic developing to the negative part of the 
invariant axis of $\rho(\gamma)$. 
 \end{ex}
 By a classical result of Goldman (see \cite[Theorem C]{GO}), given any complex 
projective structure $\sigma$ with \qf holonomy $\rho$ there exists a unique multicurve $\gamma$ on 
$S$ such that $\sigma=Gr(\sigma_\rho,\gamma)$.

\subsection{Bubbling}\label{s_bubbling}
The second surgery we will consider is a variation of the previous one, which uses a simple arc 
instead of a simple close curve, and was first considered by Gallo-Kapovich-Marden  in \cite{GKM}. 
As before, let   $\sigma$ be defined by a couple $(dev,\rho)$.
\begin{defin}
 A simple compact arc $\beta \subset S$ is bubbleable with respect to $\sigma$ if it is injectively 
developed, i.e. the restriction of $dev$ to any of its lifts $\widetilde{\beta}\subset 
\widetilde{S}$ is injective.
\end{defin}
Notice that if $\beta$ is bubbleable then the complement of its developed image in $\cp$ is a 
disk.
\begin{defin}
 Let $\beta\subset S$ be a bubbleable arc with respect to $\sigma$. For any lift 
$\widetilde{\beta}$ of $\beta$ we cut $\widetilde{S}$ along it and a copy of $\cp$ along 
$dev(\widetilde{\beta})$, and glue them together equivariantly via the developing map. 
This produces a simply connected surface $\widetilde{S}'$ to which the action 
$\pi_1(S)\curvearrowright \widetilde{S}$ and the map $dev:\widetilde{S}\to \cp$ naturally extend, 
so that the quotient is naturally endowed with a new geometric structure. We 
call this structure the \textbf{bubbling} of $\sigma$ along $\beta$ and denote it by 
$Bub(\sigma,\beta)$. The surface $\sigma \setminus \beta$  embeds in 
$Bub(\sigma,\beta)$ and the complement is the disk $B=\cp  \setminus dev(\widetilde{\beta})$, which 
we call the bubble associated to $\beta$. 
\end{defin}
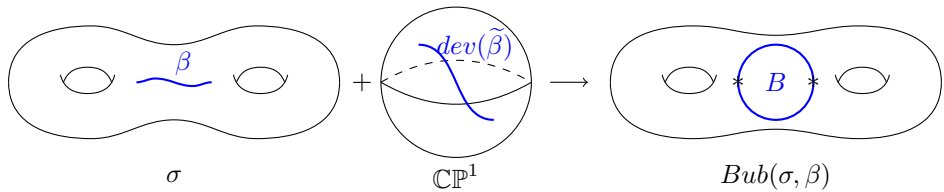
\begin{figure}[h]

\begin{center}
\begin{tikzpicture}[scale=0.5]
\node at (0,-2.5) {$\sigma$};
\draw[xscale=0.8] (-5.5,0) to[out=90,in=180] (-2.75,1.5) to[out=0,in=180] (0,1) to[out=0,in=180] (2.75,1.5) to[out=0,in=90] (5.5,0) ;
\draw[rotate=180,xscale=0.8] (-5.5,0) to[out=90,in=180] (-2.75,1.5) to[out=0,in=180] (0,1) to[out=0,in=180] (2.75,1.5) to[out=0,in=90] (5.5,0) ;

\draw[xshift=-0.3cm,yshift=-0.1cm,xscale=0.8]  (-3.3,0.1) to[out=65,in=180] (-2.5,0.5) to[out=0,in=115] (-1.7,0.1);
\draw[xshift=-0.3cm,yshift=-0.1cm,xscale=0.8]  (-3.4,0.3) to[out=-75,in=180] (-2.5,-0.2) to[out=0,in=-105] (-1.6,0.3);
\draw[xshift=4.3cm,yshift=-0.1cm,xscale=0.8]  (-3.3,0.1) to[out=65,in=180] (-2.5,0.5) to[out=0,in=115] (-1.7,0.1);
\draw[xshift=4.3cm,yshift=-0.1cm,xscale=0.8]  (-3.4,0.3) to[out=-75,in=180] (-2.5,-0.2) to[out=0,in=-105] (-1.6,0.3);

\node at (0.25,0.5) {\color{blue} $\beta$};
\draw[blue,thick] (-1,0) to[out=0,in=180] (-0.5,0.1) to[out=0,in=180] (0.5,-0.1) to[out=0,in=180] (1,0);

\node at (5,0) {$+$};

\begin{scope}[xshift=-2.5cm]
\node at (10,-2.5) {$\mathbb{CP}^1$};
\node at (10.5,1) {\color{blue} $dev(\widetilde{\beta})$};
\draw (8,0) to[out=90,in=180] (10,2) to[out=0,in=90] (12,0) to[out=-90,in=0] (10,-2)  to[out=180,in=-90] (8,0);
\draw (8,0) to[out=-30,in=210] (12,0);
\draw[dashed] (8,0) to[out=30,in=150] (12,0);
\draw[blue,thick] (9,1) to[out=0,in=180] (11,-1);
\draw[->] (12.5,0) -- (13.5,0);
\end{scope}

\begin{scope}[xshift=16cm]
\node at (0,-2.5) {$Bub(\sigma,\beta)$};
\draw[xscale=0.8] (-5.5,0) to[out=90,in=180] (-2.75,1.5) to[out=0,in=180] (0,1.25) to[out=0,in=180] (2.75,1.5) to[out=0,in=90] (5.5,0) ;
\draw[rotate=180,xscale=0.8] (-5.5,0) to[out=90,in=180] (-2.75,1.5) to[out=0,in=180] (0,1.25) to[out=0,in=180] (2.75,1.5) to[out=0,in=90] (5.5,0) ;

\draw[xshift=-0.3cm,yshift=-0.1cm,xscale=0.8]  (-3.3,0.1) to[out=65,in=180] (-2.5,0.5) to[out=0,in=115] (-1.7,0.1);
\draw[xshift=-0.3cm,yshift=-0.1cm,xscale=0.8]  (-3.4,0.3) to[out=-75,in=180] (-2.5,-0.2) to[out=0,in=-105] (-1.6,0.3);
\draw[xshift=4.3cm,yshift=-0.1cm,xscale=0.8]  (-3.3,0.1) to[out=65,in=180] (-2.5,0.5) to[out=0,in=115] (-1.7,0.1);
\draw[xshift=4.3cm,yshift=-0.1cm,xscale=0.8]  (-3.4,0.3) to[out=-75,in=180] (-2.5,-0.2) to[out=0,in=-105] (-1.6,0.3);

\draw[blue,thick] (-1,0) to[out=90,in=180] (0,1) to[out=0,in=90] (1,0);
\draw[blue,thick] (-1,0) to[out=-90,in=180] (0,-1) to[out=0,in=-90] (1,0);
\node at (-1,0) {$*$};
\node at (1,0) {$*$};
\node at (0,0) {\color{blue} $B$};

\end{scope}

\end{tikzpicture}

\end{center}
\caption{Bubbling a surface}
\end{figure}
\begin{rmk}\label{branch!}
 The reader should be warned that in the case of grafting the arc 
$\widetilde\gamma$ has endpoints at infinity (i.e. it does not have a compact 
closure in $\widetilde{S}$), whereas a bubbleable arc $\beta$ has its endpoints on 
the surface; strictly speaking a bubbling does not produce a complex projective structure: 
after it has been performed, the geometric structure branches around the endpoints of the arc. 
Such a structure is known in the literature as a branched complex projective structure. We are not 
concerned with branched structures on their own in this paper; instead we will use them just as a 
tool to study the grafting surgery on unbranched structures. Therefore we content ourselves with 
saying that the definition of this class of structures is the same as the one given in 
\ref{def_cpstructure}, but local charts are allowed to be finite orientation preserving branched 
covers; we refer to \cite{MA} (were they were first introduced) for standard background. Also 
the definition of the geometric decomposition for structures with \qf holonomy goes through as in 
the case of unbranched structures, with minor modification (see \cite{CDF} for more details).
\end{rmk}

\begin{ex}
The easiest example is obtained by bubbling a hyperbolic surface $\sigma_\rho$ along an embedded 
compact geodesic arc $\beta$. On $Bub(\sigma_\rho,\beta)$ we see a negative disk isometric to the 
lower-half plane, bounded by a simple closed curve isomorphic to $\rp$. The positive part consists 
of a subsurface isometric to $\sigma_\rho \setminus \beta$ glued along a copy of $\hyp^2 
\setminus dev(\widetilde{\beta})$; notice that this positive component contains a couple of simple 
branch points, i.e. the angle around each of them (with respect to the induced conformal structure) 
is $4\pi$.
\end{ex}
To reverse this surgery one needs to find a subsurface which can be removed, in the same way a 
grafting annulus can. We find it convenient to give the following definition.
\begin{defin}
A \textbf{bubble} on a branched complex projective structure $\sigma$ is an 
embedded closed disk $B \subset S$  whose boundary decomposes as $\partial B=\beta' \cup \{x,y\} 
\cup \beta''$ where $\{x,y\}$ are simple branch points of $\sigma$ and $\beta',\beta''$ are 
embedded 
injectively developed arcs which overlap once developed; more precisely there exist a determination 
of the developing map on $B$ which injectively maps $\beta',\beta''$ to the same simple arc 
$\widehat{\beta}\subset \cp$ and restricts to a diffeomorphism $dev: int(B)\to \cp \setminus 
\widehat{\beta} $. A \textbf{debubbling} is the surgery which consists in removing a bubble and 
gluing the resulting boundary.
\end{defin}
As in the case of grafting, both bubbling and debubbling preserve the holonomy of the structure.
A result analogous to Goldman's theorem for branched structure was obtained by the author (see 
\cite[Theorem 1.1]{R}), namely that a generic branched complex projective structure with \qf 
holonomy and two simple branch points is obtained by bubbling an unbranched structure with the same 
holonomy.

\section{Degrafting via bubbles}\label{s_bubgraft}
A relation between the two surgeries introduced above has first been obtained in \cite[Theorem 
5.1]{CDF}, where it is shown that the grafting of a complex projective structure $\sigma$ 
along a simple graftable curve $\gamma$ can always be obtained by performing first a bubbling and 
then a debubbling; more precisely if the first bubbling is performed along a bubbleable arc 
$\beta\subset \sigma$, then $Bub(\sigma,\beta)$ displays a bubble $B$ coming from $\beta$, and the 
content of the theorem is that it is possible to find a different bubble $B'$ corresponding to some 
bubbleable arc $\beta'\subset Gr(\sigma, \gamma)$ and such that 
$Bub(\sigma,\beta)=Bub(Gr(\sigma,\gamma),\beta')$. The two bubbles $B,B'$ are not 
isotopic relative to the branch points. Here we prove that in \qf 
holonomy this procedure also works for a general multi(de)grafting, i.e. that grafting a multicurve 
(or degrafting a collection of grafting annuli) can always be realised by a simple sequence of one 
bubbling and one debubbling.\par
To simplify the exposition we adopt the convention that normal letters denote objects on the 
surface, letters with a tilde denote a lift to the universal cover and letters with a hat 
denote a developed image of the corresponding object.\par
Let us fix a Fuchsian representation $\rho : \pi_1(S)\to \pslr$ and a projective structure $\sigma$ 
with holonomy $\rho$. By \cite{GO} $\sigma$ is obtained by a multigrafting on the  uniformizing 
structure $\sigma_\rho=\hyp^2/\rho(\pi_1(S))$, hence it decomposes into a 
hyperbolic core of finite volume (coming directly from $\sigma_\rho$) plus a certain number of 
grafting regions. We will denote by $A_\gamma=A^1_\gamma \cup \dots \cup A^M_\gamma$  the grafting 
region obtained by grafting $M$ times some simple closed geodesic $\gamma$ of $\sigma_\rho$.   
Recall from \ref{ex_graftgeod} that the geometric decomposition of a grafting annulus is made of a 
negative annulus and a couple of ends in the adjacent positive component(s). Notice that the 
structure on the interior of each grafting annulus is uniformizable, in the sense that the 
developing image is injective on the interior of the universal cover of the annulus.\par

\subsection{Crossing a grafting annulus}\label{s_annulus}
In this section we show how to remove a simple grafting annulus on a projective structure by 
bubbling and debubbling it, by proving a more general statement. In the previous notation, let us 
begin by considering the case $M=1$ and let $A_\gamma=A_\gamma^1$ be a grafting annulus, with 
boundary geodesics $\gamma_L,\gamma_R$. Given a bubbleable arc $\beta$ which crosses $A_\gamma$ 
transversely from side to side, we introduce some auxiliary objects which are needed in the main 
construction.  

\begin{defin}
 Let $\beta$ be an oriented bubbleable arc properly embedded in $A_\gamma$ (i.e. $\partial \beta = 
\beta \cap \partial A_\gamma$). We  call $I$ (in) and $O$ (out) the two points of $\partial 
\beta $ at which $\beta$ respectively enters in the annulus and leaves it.
Notice that there is a unique point on $\partial A_\gamma$ which is different from $I$ but is 
developed to the same point $\widehat{I}$. We will refer to it as the twin of $I$, and 
similarly for $O$.
\end{defin}
\begin{defin}
Let $\beta$ be an oriented bubbleable arc properly embedded in $A_\gamma$. We  define a 
preferred orientation for $\gamma$ so that in the developed image $\widehat{O}$ sits after 
$\widehat{I}$ along $\widehat{\gamma}$ (since $\beta$ is bubbleable, $\widehat{I}\neq \widehat{O}$, 
thus this is well defined). We refer to it as the orientation of $\gamma$ induced by $\beta$.
\end{defin}
\begin{figure}[h]

\begin{center}
\begin{tikzpicture}

\draw (-4.5,2) -- (-1.5,2) -- (-1.5,-2) -- (-4.5,-2) -- (-4.5,2);
\draw (-5,2) -- (-4,2);
\draw (-1,2) -- (-2,2);
\draw (-5,-2) -- (-4,-2);
\draw (-1,-2) -- (-2,-2);
\node at (-3,1) {\color{blue} $\beta$};
\draw[blue,thick,->] (-4.5,-1) -- (-3,0);
\draw[blue,thick] (-3,0) -- (-1.5,1);
\node at (-5,-1) { $I$};
\node at (-1,-1) { $I^{twin}$};
\node at (-1,1) { $O$};
\node at (-5.2,1) { $O^{twin}$};
\node at (-4.5,-1) { $\bullet$};
\node at (-1.5,-1) { $\bullet$};
\node at (-1.5,1) { $\bullet$};
\node at (-4.5,1) { $\bullet$};
\node at (-3,-3) { $A_\gamma$};
\node at (-4.5,-2.5) { $\gamma_L$};
\node at (-1.5,-2.5) { $\gamma_R$};

 \draw (1,0) -- (5,0);
\draw[->] (3,-2) -- (3,2);
\node at (3.5,-2) { $\widehat{\gamma}$};
\node at (4.5,-1) {\color{blue} $\widehat{\beta}$};
\node at (2.5,0.5) { $\widehat{I}$};
\node at (3.5,1.5) { $\widehat{O}$};
\node at (3,0.5) { $\bullet$};
\node at (3,1.5) { $\bullet$};
\node at (3,-3) { $\cp$};

\draw[blue,thick,->] (3,0.5) to[out=0,in=90] (3.75,0) to[out=-90,in=0] (3,-1);
\draw[blue,thick] (3,-1) to[out=180,in=-90] (2,0)  to[out=90,in=180] (3,1.5);

\end{tikzpicture}
\end{center}
\caption{An arc inducing In and Out points and an orientation.}
\end{figure}
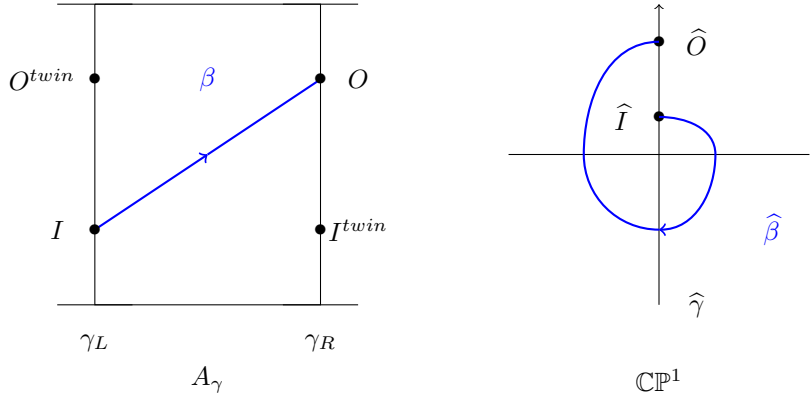

\begin{defin}
 Let $\beta, \beta'$ be oriented bubbleable arcs properly embedded in $A_\gamma$. We say that 
$\beta'$ is \textbf{coherent} with $\beta$ if $\widehat{O'}$ sits after $\widehat{I'}$ along 
$\widehat{\gamma}$ with respect to the orientation induced by $\beta$. Otherwise we say that it 
is incoherent. 
\end{defin}
\begin{figure}[h]

\begin{center}
\begin{tikzpicture}

\draw (-4.5,2) -- (-1.5,2) -- (-1.5,-2) -- (-4.5,-2) -- (-4.5,2);
\draw (-5,2) -- (-4,2);
\draw (-1,2) -- (-2,2);
\draw (-5,-2) -- (-4,-2);
\draw (-1,-2) -- (-2,-2);
\draw[blue,thick,->] (-4.5,-1.5) -- (-3,-1);
\draw[blue,thick] (-3,-1) -- (-1.5,-0.5);
\draw[blue,thick,->] (-4.5,0.5) -- (-3,1);
\draw[blue,thick] (-3,1) -- (-1.5,1.5);
\node at (-3,-3) { $A_\gamma$};
\node at (-4.5,-2.5) { $\gamma_L$};
\node at (-1.5,-2.5) { $\gamma_R$};
\node at (-3,1.5) {\color{blue} $\beta'$};
\node at (-3,-1.5) {\color{blue} $\beta$};

\begin{scope}[xshift=6cm]
\draw (-4.5,2) -- (-1.5,2) -- (-1.5,-2) -- (-4.5,-2) -- (-4.5,2);
\draw (-5,2) -- (-4,2);
\draw (-1,2) -- (-2,2);
\draw (-5,-2) -- (-4,-2);
\draw (-1,-2) -- (-2,-2);
\draw[blue,thick,->] (-4.5,-1.5) -- (-3,-1);
\draw[blue,thick] (-3,-1) -- (-1.5,-0.5);
\draw[blue,thick,->] (-4.5,1.5) -- (-3,1);
\draw[blue,thick] (-3,1) -- (-1.5,0.5);
\node at (-3,-3) { $A_\gamma$};
\node at (-4.5,-2.5) { $\gamma_L$};
\node at (-1.5,-2.5) { $\gamma_R$};
\node at (-3,1.5) {\color{blue} $\beta'$};
\node at (-3,-1.5) {\color{blue} $\beta$};
\end{scope}

\end{tikzpicture}
\end{center}

\caption{Coherent and incoherent arcs.}
\end{figure}
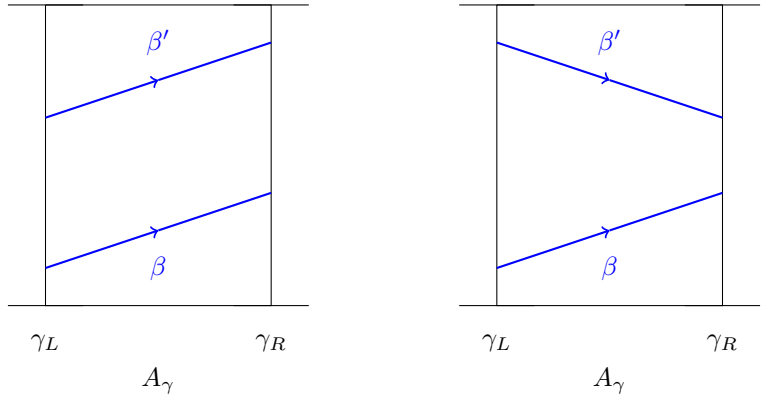
Now let $\beta$ be an oriented bubbleable arc which transversely crosses some grafting annulus 
$A_\gamma$, i.e. every time it enters $A_\gamma$ it crosses it and leaves it on the other side. 
Then $\beta \cap A_\gamma=\beta_1 \cup \dots \cup \beta_N $ is a disjoint union of oriented 
bubbleable arcs properly embedded in $A_\gamma$, which we will call crossings; the labelling of 
$\beta_1,\dots,\beta_N$ is such that they appear in this order along $\beta$. For each  
 crossing $\beta_k$ we can define the entry and exit point $I_k$ and $O_k$, and the induced 
orientation of $\gamma$ as above. Since $\beta$ is 
embedded and bubbleable all these points are distinct, and the same holds for their developed 
images. We agree to fix the orientation of $\gamma$ determined by the first crossing $\beta_1$, but 
of course more generally we can also decide if two given crossings are coherent or not with respect 
to each other. \par
Let us introduce another useful way to order the crossings $\beta_1,\dots,\beta_N$, 
according to the way they appear when travelling along $\gamma$ with respect to the orientation of 
$\gamma$ induced by $\beta_1$: set $\alpha_1=\beta_1$, then let $\alpha_{k+1}$ be the crossing 
we meet after $\alpha_k$ along $\gamma$ with respect to the chosen orientation. We get an ordering 
of the crossings as $\alpha_1,\dots,\alpha_N$ which is actually a $\intz_N$-order (i.e. 
$\alpha_{N+1}=\alpha_1$); moreover there exists a unique permutation $\sigma \in \mathfrak{S}_N$ 
such that $\alpha_k=\beta_{\sigma(k)}$ and $\sigma(1)=1$. 
We keep track of the coherence between crossings by defining the following \textbf{coherence 
parameters}
$$\varepsilon_k=\left\lbrace \begin{array}{ll}
                              1 & \textrm{if $\alpha_k$ coherent with $\alpha_1=\beta_1$}\\
                              -1 & \textrm{if $\alpha_k$ incoherent with $\alpha_1=\beta_1$}\\
                             \end{array}\right. $$
$$\varepsilon_{k,l}=\left\lbrace \begin{array}{ll}
                              1 & \textrm{if $\alpha_k,\alpha_l$ coherent with each other}\\
                              -1 & \textrm{if $\alpha_k,\alpha_l$ incoherent with each other}\\
                             \end{array}\right. $$

Let us roughly describe the idea behind the main construction of this section. 
Given a bubbleable arc which transversely crosses a grafting annulus, we would like to perform the 
bubbling along it and then find another bubble which avoids the real curve. The naive approach is 
to 
start from a branch point and follow the given bubble until we meet the region corresponding to the 
grafting annulus at the points coming from $I_1$; here one path can follow the curve coming from 
the 
boundary of the grafting annulus until the twin of $O_1$, and the other one can follow its analytic 
extension inside the bubble to cross the bubble from side to side. Notice that in doing this  it 
also crosses the grafting annulus from side to side; in particular it reaches $O_1$. Then they keep 
travelling along the boundaries of the grafting annulus in the direction induced by $\beta_1$, 
until 
they meet $\alpha_2$. One of them will meet that crossing before the other and will follow the 
analytic extension of $\gamma$ inside the bubble, while the other one will follow the boundary of 
the annulus; the coherence parameters $\varepsilon_{k}$ and $\varepsilon_{k,l}$ determine the order 
in which points are met, and the direction in 
which the paths will go. Anyway they will reach points on the same side of the annulus, but on 
opposite sides of the bubble, hence they can keep walking along the original bubble. This works 
because at every crossing there is an analytic extensions of $\gamma$ inside the bubble which 
crosses it from side to side. However 
in general this naive procedure does not result in a couple of disjoint embedded arcs: already in 
the case of a single crossing ($N=1$) the analytic extension of $\gamma$ inside the bubble is used 
twice; more precisely, the two arcs we naively construct in this way intersect along the analytic 
extension of $\gamma$ inside the bubble, hence we do not get a new bubble, but only an embedded 
open disk whose boundary self-bumps along an arc.\par

To fix this we consider a small collar neighbourhood $A_\gamma^{\#}$ of $A_\gamma$; this can be 
obtained by slightly pushing the boundary curves of $A_\gamma$ into the hyperbolic core of the 
adjacent components (i.e. away from $A_\gamma$). More precisely it can be taken to be the region 
bounded by the couple of simple closed curves $\gamma_{\pm 1} = \{ x \in S^+\setminus A_\gamma \ | 
\ d(x,\gamma)= \varepsilon \}$, for some small $\varepsilon >0$, which develop to the two 
boundaries of the region $\mathcal{N}_\varepsilon (\widehat{\gamma})=\{ \widehat{x} \in \hyp^2 \ | 
\ d(\widehat{x},\widehat{\gamma})\leq \varepsilon \}$. Notice that 
the developing image is no longer injective in the interior of $A_\gamma^\#$.\par

We have that $\beta \cap A_\gamma^\#=\beta_1^\# \cup \dots \cup \beta_N^\# $ is a disjoint 
union of oriented bubbleable arcs properly embedded in $A_\gamma^\#$ and such that 
$\beta_k\subset \beta^\#_k$, which we still call crossing. Moreover each crossing $\beta_k^\#$ will 
intersect $\partial A_\gamma^{\#}$ in two points; let us label them by  $I_k^{-1}$ and 
$O_k^{+1}$ in such a way that  $I_k^{-1},I_k,O_k,O_k^{+1}$ appear in this order  along $\beta$; 
then 
label the curves $\gamma_{\pm 1}$ so that $I_1^{-1}\in \gamma_{-1}$ and $O_1^{+1}\in \gamma_{+1}$. 
Notice that for the other crossings it may happen that $I_k^{-1}\in \gamma_{\pm1}$ and $O_k^{+1}\in 
\gamma_{\mp1}$, according to the fact that $\beta^\#_k$ enters in $A_\gamma^\#$ on the same side as
$\beta_1^\#$ leaves it or not; however this is not going to be a relevant in our construction.\par

 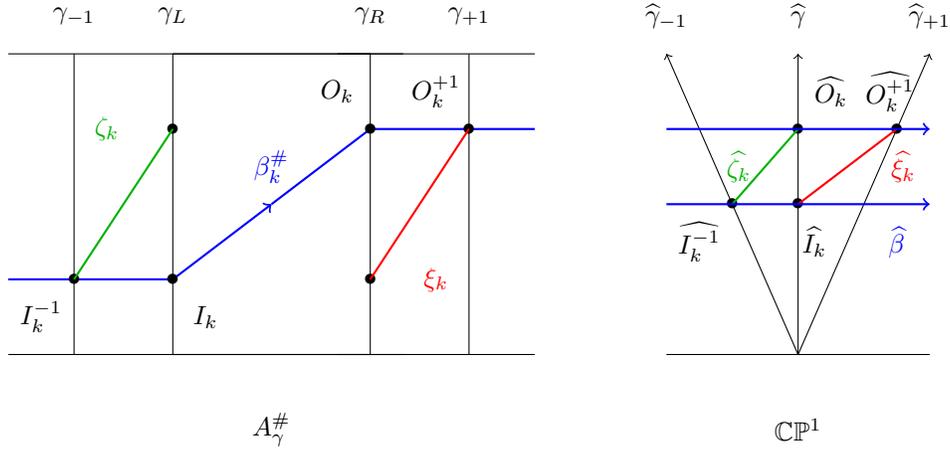
\begin{figure}[h]

\begin{center}
\begin{tikzpicture}[xscale=0.875]

\draw (-4.5,2) -- (-1.5,2) -- (-1.5,-2) -- (-4.5,-2) -- (-4.5,2);
\draw (-7,2) -- (1,2);
\draw (-1,2) -- (-2,2);
\draw (-7,-2) -- (1,-2);
\draw (-1,-2) -- (-2,-2);
\draw (-6,-2) -- (-6,2);
\draw (0,-2) -- (0,2);
\node at (-3,0.5) {\color{blue} $\beta_k^\#$};
\draw[blue,thick,->] (-7,-1) to[out=0,in=180] (-4.5,-1) -- (-3,0);
\draw[blue,thick] (-3,0) -- (-1.5,1) to[out=0,in=180] (1,1);
\node at (-4,-1.5) { $I_k$};
\node at (-6.5,-1.5) { $I_k^{-1}$};
\node at (-2,1.5) { $O_k$};
\node at (-0.5,1.5) { $O_k^{+1}$};
\node at (-4.5,-1) { $\bullet$};
\node at (-6,-1) { $\bullet$};
\node at (-1.5,-1) { $\bullet$};
\node at (0,1) { $\bullet$};
\node at (-1.5,1) { $\bullet$};
\node at (-4.5,1) { $\bullet$};
\node at (-3,-3) { $A_\gamma^\#$};
\node at (-6,2.5) { $\gamma_{-1}$};
\node at (0,2.5) { $\gamma_{+1}$};
\node at (-4.5,2.5) { $\gamma_L$};
\node at (-1.5,2.5) { $\gamma_R$};

\draw[black!30!green,thick] (-6,-1) -- (-4.5,1);
\draw[red,thick] (-1.5,-1) -- (0,1);
\node at (-5.5,0.2) {\color{black!30!green} $\zeta_k$};
\node at (-0.5,-0.2) {\color{red} $\xi_k$};

 \draw (3,-2) -- (7,-2);
\draw[->] (5,-2) -- (5,2);
\draw[->] (5,-2) -- (3,2);
\draw[->] (5,-2) -- (7,2);
\node at (5,2.5) { $\widehat{\gamma}$};
\node at (3,2.5) { $\widehat{\gamma}_{-1}$};
\node at (7,2.5) { $\widehat{\gamma}_{+1}$};

\node at (5,0) { $\bullet$};
\node at (4,0) { $\bullet$};
\node at (5,1) { $\bullet$};
\node at (6.5,1) { $\bullet$};

\node at (5.25,-0.5) { $\widehat{I_k}$};
\node at (3.5,-0.5) { $\widehat{I_k^{-1}}$};
\node at (5.5,1.5) { $\widehat{O_k}$};
\node at (6.4,1.5) { $\widehat{O_k^{+1}}$};

\node at (4.1,0.5) {\color{black!30!green} $\widehat{\zeta_k}$};
\node at (6.6,0.5) {\color{red} $\widehat{\xi_k}$};

\draw[black!30!green,thick] (4,0) -- (5,1);
\draw[red,thick] (5,0) -- (6.5,1);

\node at (6.5,-0.5) {\color{blue} $\widehat{\beta}$};

\node at (5,-3) { $\cp$};

\draw[blue,thick,->] (3,0) to[out=0,in=180] (7,0) ;
\draw[blue,thick,->] (3,1) to[out=0,in=180] (7,1) ;

\end{tikzpicture}
\end{center}
 \caption{The extended annulus $A_\gamma^\#$ and the auxiliary objects.}
 \end{figure}

Now for any $k=1,\dots, N$ we consider in the developed image in $\hyp^2$ the geodesic segment 
$\widehat{\zeta_k}$ from $\widehat{I_k^{-1}}$ to $\widehat{O_k}$ and the geodesic segment 
$\widehat{\xi_k}$ from $\widehat{I_k}$ to $\widehat{O_k^{+1}}$. This defines for us an arc 
$\zeta_k$  in $A_\gamma^{\#} \setminus A_\gamma$ starting from $I_k^{-1}$ and ending at the twin of 
 $O_k$, and an arc $\xi_k$ in $A_\gamma^{\#} \setminus A_\gamma$ starting from the twin of  $I_k$ 
 and ending at $O_k^{+1}$. Since $\beta$ is embedded and bubbleable, all these arcs are disjoint; 
notice that the behaviour of $\zeta_k$ and $\xi_k$ in $A_\gamma^\#$ essentially mimics that of 
$\beta_k$ (e.g. they wrap around $A_\gamma$ the same number of times), with the only difference 
that they are entirely contained in the positive region, while 
$\beta_k$ crosses the real curve twice inside $A_\gamma$.
To simplify the exposition we also find it convenient to introduce an action of $\intz_2=\{\pm 1\}$ 
on all the auxiliary objects we have defined: we let $1$ act as the identity, while $-1$ acts by 
exchanging an ``entry object'' with the corresponding ``exit object'', i.e.
$$-1.I_k=O_k \quad -1.I_k^{-1}=O_k^{+1} \quad -1.\zeta_k=\xi_k$$
Moreover notice that all arcs involved are oriented; for any path 
$\mu$, let $\mu^{-1}$ denote the same path with the opposite orientation.
We now have all the ingredients required to prove the following result. 
\begin{prop}\label{onegraftingannulus}
Let $\rho:\pi_1(S)\to \pslr$ be Fuchsian.
 Let $\sigma_0$ be a complex projective structure with holonomy $\rho$, $A_\gamma \subset \sigma_0$ 
a grafting annulus and $\beta \subset \sigma_0$ a bubbleable arc which transversely crosses 
$A_\gamma$ and avoids all other grafting annuli of $\sigma_0$. Then there exist a complex 
projective structure $\sigma_0'$ with the same holonomy $\rho$ and a bubbleable arc $\beta'\subset 
\sigma_0'$  which avoids all the real curves of $\sigma_0'$ and such that 
$Bub(\sigma_0,\beta)=Bub(\sigma_0',\beta')$.
\end{prop}
\proof
We will prove this by directly finding a new bubble with the required properties on 
$Bub(\sigma_0,\beta)$. Pick an orientation of $\beta$; then we have all the auxiliary objects 
defined above, in particular fix the orientation of $\gamma$ induced by the first crossing 
$\alpha_1=\beta_1$. We will define a new bubble roughly in the following way: 
each time the bubble coming from $\beta$ enters $A_\gamma^\#$ in correspondence of some crossing 
$\alpha_k$ we will describe how to leave $A_\gamma^\#$ in correspondence of the crossing 
$\alpha_{k\pm 1}$ by suitably following some of the auxiliary arcs (the sign depends on some 
coherence parameters); then we keep following $\beta$ until we reach another crossing, if any, and 
we iterate. \par
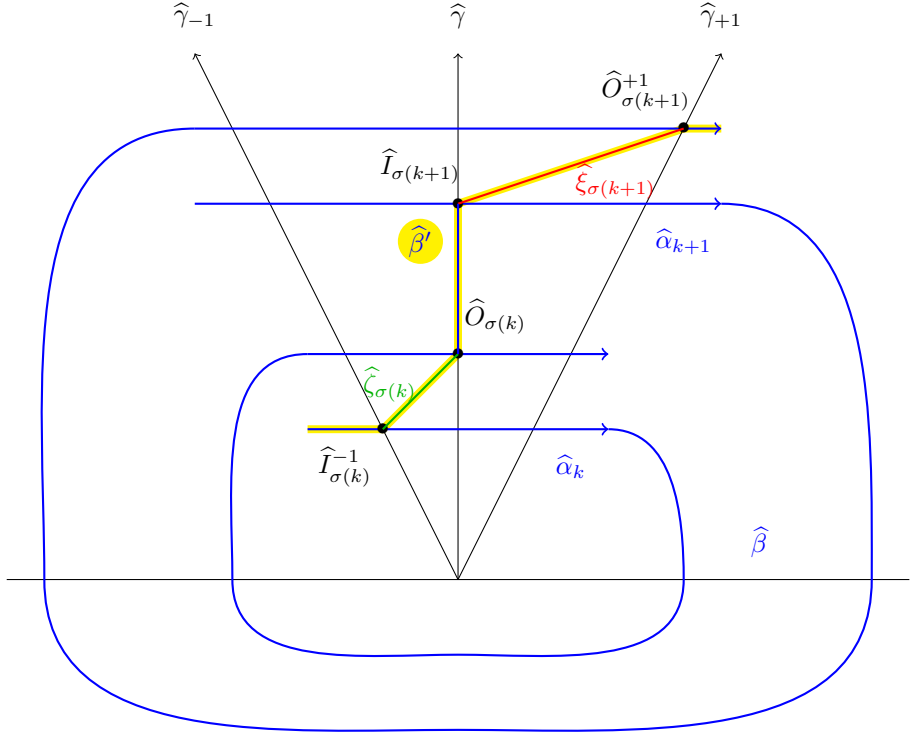
\begin{figure}[h]

\begin{center}
\begin{tikzpicture}[scale=1]

\draw[yellow,line width=3pt] (3,0) -- (4,0) -- (5,1) -- (5,3) -- (8,4) -- (8.5,4);
\fill[yellow] (4.5,2.5) circle (.3cm);

 \draw (-1,-2) -- (11,-2);
\draw[->] (5,-2) -- (5,5); \node at (5,5.5) { $\widehat{\gamma}$};
\draw[->] (5,-2) -- (1.5,5); \node at (1.5,5.5) { $\widehat{\gamma}_{-1}$};
\draw[->] (5,-2) -- (8.5,5); \node at (8.5,5.5) { $\widehat{\gamma}_{+1}$};

\draw[blue,thick,->] (3,0) to[out=0,in=180] (7,0) ;
\draw[blue,thick,->] (3,1) to[out=0,in=180] (7,1) ;
\node at (6.5,-0.5) {\color{blue} $\widehat{\alpha}_k$};

\node at (4,0) { $\bullet$};\node at (3.5,-0.5) { $\widehat{I}_{\sigma(k)}^{-1}$};
\node at (5,1) { $\bullet$};\node at (5.5,1.5) { $\widehat{O}_{\sigma(k)}$};

\draw[black!30!green,thick] (4,0) -- (5,1); 
\node at (4.1,0.6) {\color{black!30!green} $\widehat{\zeta}_{\sigma(k)}$};


\draw[blue,thick,->] (1.5,3) to[out=0,in=180] (8.5,3) ;
\draw[blue,thick,->] (1.5,4) to[out=0,in=180] (8.5,4) ;
\node at (8,2.5) {\color{blue} $\widehat{\alpha}_{k+1}$};

\node at (5,3) { $\bullet$};\node at (4.5,3.5) { $\widehat{I}_{\sigma(k+1)}$};
\node at (8,4) { $\bullet$};\node at (7.5,4.5) { $\widehat{O}_{\sigma(k+1)}^{+1}$};

%
\draw[red,thick] (5,3) -- (8,4);
\node at (7.1,3.3) {\color{red} $\widehat{\xi}_{\sigma(k+1)}$};

\draw[blue,thick] (5,1) -- (5,3); \node at (4.5,2.5) {\color{blue} $\widehat{\beta}'$};

\draw[blue,thick] (7,0) to[out=0,in=90] (8,-2) to[out=-90,in=0] (5,-3) to[out=180,in=-90] (2,-2) 
to[out=90,in=180] (3,1);

\draw[blue,thick] (8.5,3) to[out=0,in=90] (10.5,-2) to[out=-90,in=0] (5,-4) to[out=180,in=-90] 
(-.5,-2) to[out=90,in=180] (1.5,4);

\node at (9,-1.5) {\color{blue} $\widehat{\beta}$};

\end{tikzpicture}
\end{center}
  \caption{The path $\widehat{\beta}'$ in $\cp$: the $k$-th crossing in the case 
 $\omega = \varepsilon_k = \varepsilon_{k,k+1} = 1$ and $\widehat{I}_{\sigma(k+1)}$ sits after 
 $\widehat{O}_{\sigma(k)}$ along $\widehat{\gamma}$.}\label{pic_devpath}
 \end{figure}

Let us now define a procedure to handle the $k$-th crossing in the developed image (see Picture 
\ref{pic_devpath}). 
Suppose $\widehat{\beta}$ enters  $\mathcal{N}_\varepsilon(\widehat{\gamma})$ in correspondence of 
$\widehat{\alpha}_k=\widehat{\beta}_{\sigma(k)}$  at 
$\omega \widehat{I}_{\sigma(k)}^{-1}$ for some $\omega \in \{\pm 1\}$. 
We begin by following $\omega \widehat{\zeta}_{\sigma(k)}^\omega$, so that we get to $\omega 
\widehat{O}_{\sigma(k)}$. We now distinguish two cases according to the relative position of the 
endpoints of the two crossings $\widehat{\alpha}_k$ and $\widehat{\alpha}_{k+\omega \varepsilon_k}$
\begin{enumerate}
 \item if  $\omega \varepsilon_{k,k+\omega 
\varepsilon_k}\widehat{I}_{\sigma(k+\omega\varepsilon_k)}$ sits after $\omega 
\widehat{O}_{\sigma(k)}$ along $\widehat{\gamma}^{\omega \varepsilon_k}$, then we follow 
$\widehat{\gamma}^{\omega \varepsilon_k}$ until we reach it; we meet 
$\widehat{\alpha}_{k+\omega \varepsilon_k}$ at that point $\omega \varepsilon_{k,k+\omega 
\varepsilon_k}\widehat{I}_{\sigma(k+\omega\varepsilon_k)}$ and then we can follow the arc $\omega 
\varepsilon_{k,k+\omega \varepsilon_k}\widehat{\xi}_{\sigma(k+\omega\varepsilon_k)}^{\omega 
\varepsilon_{k,k+\omega \varepsilon_k}}$
\item otherwise  $\omega \varepsilon_{k,k+\omega 
\varepsilon_k}\widehat{I}_{\sigma(k+\omega\varepsilon_k)}$ sits 
before $\omega \widehat{O}_{\sigma(k)}$ along $\widehat{\gamma}^{\omega \varepsilon_k}$, then the 
fact that $\beta$ is embedded implies that $-\omega \varepsilon_{k,k+\omega 
\varepsilon_k}\widehat{I}_{\sigma(k+\omega\varepsilon_k)}$ is 
after $\omega \widehat{O}_{\sigma(k)}$; in this case we can move a little off $\widehat{\gamma}$ 
along $\widehat{\beta}^\omega$ to meet the arc
$\omega \varepsilon_{k,k+\omega \varepsilon_k}\widehat{\xi}_{\sigma(k+\omega\varepsilon_k)}^{\omega 
\varepsilon_{k,k+\omega \varepsilon_k}}$
\end{enumerate}
In both cases we follow the arc $\omega \varepsilon_{k,k+\omega 
\varepsilon_k}\widehat{\xi}_{\sigma(k+\omega\varepsilon_k)}^{\omega \varepsilon_{k,k+\omega 
\varepsilon_k}}$ and reach $\omega \varepsilon_{k,k+\omega 
\varepsilon_k}\widehat{O}_{\sigma(k+\omega\varepsilon_k)}^{+1}$. Then we are ready to leave 
$\mathcal{N}_\varepsilon(\widehat{\gamma})$ along $\beta^{\omega \varepsilon_{k,k+\omega 
\varepsilon_k}}$.  We use this rule to define a path $\widehat{\beta}'$ in $\cp$,  starting from 
the 
first endpoint of $\widehat{\beta}$.\par
We should explicitly remark that it is possible that $\beta$ goes around some topology of the 
surface between two crossings $\beta_k$ and $\beta_{k+1}$; in this case its developed image does 
not come back to the region $\mathcal{N}_\varepsilon (\widehat{\gamma})$, but to a different region 
$g\mathcal{N}_\varepsilon (\widehat{\gamma})$ for some M\"obius transformation which depends on the 
topology around which $\beta$ travels between  $\beta_k$ and $\beta_{k+1}$. However  
translating $\mathcal{N}_\varepsilon (\widehat{\gamma})$  with the holonomy of the structure 
does not produce overlaps; this follows from the fact that the developed images of the geodesic 
$\gamma$ for the underlying uniformizing structure $\sigma_\rho$  are disjoint and the 
fact that $\varepsilon >0$ can be chosen to be arbitrarily small.
On the other hand, if the path does not go around topology (so that $\widehat{\beta}$ keeps 
intersecting the same region $\mathcal{N}_\varepsilon (\widehat{\gamma})$), then it is enough to 
notice that the above procedure is completely reversible, in the sense that at any point the 
knowledge of what arc we have used at the most recent step is enough to know what arc to use to 
perform the next one, and vice versa. This implies that the path $\widehat{\beta}'$ which is 
constructed by the above rules does not pass more than once through any of its points.\par
Finally let us consider what happens to the parts of $\widehat{\beta}'$ which are outside the 
region $\mathcal{N}_\varepsilon (\widehat{\gamma})$ and its translates. By construction they come 
from portions of $\beta$ which are outside the grafting annulus $A_\gamma$; moreover by hypothesis 
$\beta$ does not intersect other grafting annuli. Therefore the developed images of these arcs are 
the same they would be in the underlying uniformizing structure $\sigma_\rho$, in particular they 
are all disjoint. This proves that the path $\widehat{\beta}'$ is 
embedded in $\cp$. Moreover since the number of marked points ( $\omega I^{-1}_k$ and $\omega I_k$ 
) is finite, it definitively reaches the point $\widehat{O}_N^{+1}$. After that point we 
keep following $\widehat{\beta}$ till the end, i.e. its second endpoint. To sum up, 
$\widehat{\beta}'$ is an embedded path with the same endpoints as $\widehat{\beta}$ but entirely 
contained in $\hyp^2$. \par
\begin{figure}[h]
\usetikzlibrary{patterns}
\begin{center}
\begin{tikzpicture}

\fill[pattern color = gray!70, pattern = north west lines, opacity=0.8] (-4,2) -- (-1.5,2) -- 
(1.5,4.5) -- (4,4.5) -- (4,3) -- (1.5,3) -- (-1.5,.5) -- (-4,0.5) -- (-4,2);
\fill[yshift=-5cm,pattern color = gray!70, pattern = north west lines, opacity=0.8] (-4,2) 
-- (-1.5,2) -- (1.5,4.5) -- (4,4.5) -- (4,3) -- (1.5,3) -- (-1.5,.5) -- (-4,0.5) -- (-4,2);

\fill[pattern color = gray!70, pattern = north east lines, opacity=0.8] (-4,-4.5) -- (-3,-4.5) .. 
controls (0,-3) .. (1.5,-.5) -- (1.5,0) -- (-1.5,0) --  (-1.5,-.5) -- (-3,-3) -- (-4,-3) 
-- (-4,-4.5);
\fill[xscale=-1,yscale=-1,pattern color = gray!70, pattern = north east lines, opacity=0.8] 
(-4,-4.5) -- (-3,-4.5) .. 
controls (0,-3) .. (1.5,-.5) -- (1.5,0) -- (-1.5,0) --  (-1.5,-.5) -- (-3,-3) -- (-4,-3) 
-- (-4,-4.5);


\draw (-4,5) -- (4,5);
\draw (-4,2) -- (-1.5,2) -- (1.5,4.5) -- (4,4.5);
\draw (-4,.5) -- (-1.5,.5) -- (1.5,3) -- (4,3);
\draw (-4,-3) -- (-1.5,-3) -- (1.5,-.5) -- (4,-.5);
\draw (-4,-4.5) -- (-1.5,-4.5) -- (1.5,-2) -- (4,-2);
\draw (-4,-5) -- (4,-5);

\draw (-3,5) -- (-3,2);\draw (-3,.5) -- (-3,-3);\draw (-3,-4.5) -- (-3,-5);
\begin{scope}[xshift=+1.5cm]
\draw (-3,5) -- (-3,2);\draw (-3,.5) -- (-3,-3);\draw (-3,-4.5) -- (-3,-5);
\end{scope}
\begin{scope}[yscale=-1,xshift=+4.5cm]
\draw (-3,5) -- (-3,2);\draw (-3,.5) -- (-3,-3);\draw (-3,-4.5) -- (-3,-5);
\end{scope}
\begin{scope}[yscale=-1,xshift=+6cm]
\draw (-3,5) -- (-3,2);\draw (-3,.5) -- (-3,-3);\draw (-3,-4.5) -- (-3,-5);
\end{scope}

\draw[black!30!green,thick] (-3,-3) -- (-1.5,-.5);
\draw[black!30!green,thick] (-3,-4.5) .. controls (0,-3) .. (1.5,-.5);
\begin{scope}[xscale=-1,yscale=-1]
\draw[red,thick] (-3,-3) -- (-1.5,-.5);
\draw[red,thick] (-3,-4.5) .. controls (0,-3) .. (1.5,-.5);
\end{scope}

\node at (-5,-3.75) { $\alpha_k$};
\node at (-5,1.25) { $\alpha_{k+1}$};
\node at (-3,-4.5) { $\bullet$};\node at (-3.5,-4) { $I^{-1}_{\sigma(k)}$};
\node at (-3,-3) { $\bullet$};\node at (-3.5,-2.5) { $I^{-1}_{\sigma(k)}$};
\node at (-1.5,-.5) { $\bullet$};\node at (-2,-.5) { $O^{twin}_{\sigma(k)}$};
\node at (1.5,-.5) { $\bullet$};\node at (2,-1) { $O_{\sigma(k)}$};
\node at (-2.5,-1.5) {\color{black!30!green} $\zeta_{\sigma(k)} $};
\node at (.5,-3.5) {\color{black!30!green} $\zeta_{\sigma(k)} $};
\node at (0,-6) { $A_\gamma^\#$};
\node at (-1.5,-5.5) { $\gamma_L$};
\node at (1.5,-5.5) { $\gamma_R$};
\node at (-3,-5.5) { $\gamma_{-1}$};
\node at (3,-5.5) { $\gamma_{+1}$};

\begin{scope}[xscale=-1,yscale=-1]
\node at (-3,-4.5) { $\bullet$};\node at (-3.5,-4) { $O^{+1}_{\sigma(k+1)}$};
\node at (-3,-3) { $\bullet$};\node at (-3.6,-2.5) { $O^{+1}_{\sigma(k+1)}$};
\node at (-1.5,-.5) { $\bullet$};\node at (-2.5,-.5) { $I^{twin}_{\sigma(k+1)}$};
\node at (1.5,-.5) { $\bullet$};\node at (2,-1) { $I_{\sigma(k+1)}$};
\node at (.5,-3.5) {\color{red} $\xi_{\sigma(k+1)} $};
\end{scope}
\node at (3,1.5) {\color{red} $\xi_{\sigma(k+1)} $};

\end{tikzpicture}
\end{center}
 \caption{The new bubble on the surface: the $k$-th crossing in the case 
 $\omega = \varepsilon_k = \varepsilon_{k,k+1} = 1$ and $\widehat{I}_{\sigma(k+1)}$ sits after 
 $\widehat{O}_{\sigma(k)}$ along $\widehat{\gamma}$. (The two bubbles are shaded at different 
angles).}\label{pic_newbubble}
 \end{figure}
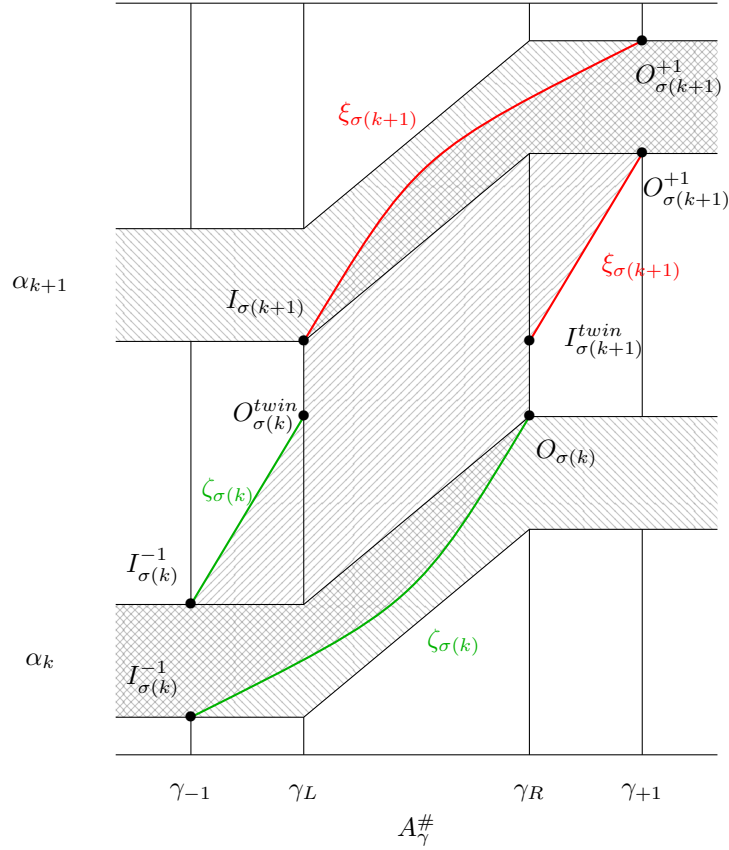
We can now follow this path on $Bub(\sigma_0,\beta)$ to identify a new bubble (see Figure 
\ref{pic_newbubble}). We start at the 
branch point of $Bub(\sigma_0,\beta)$ which is the first with respect to the chosen orientation of 
$\beta$ and follow the two twin paths developing to $\widehat{\beta}$ which give the boundary of 
the natural bubble of $Bub(\sigma_0,\beta)$. Then we check that at each crossing 
$\alpha_k$ there is a couple of embedded arcs developing to subarcs of $\widehat{\beta}'$ which 
fellow travel from the entry point $\omega I_{\sigma(k)}^{-1}$ to the exit point $\omega 
\varepsilon_{k,k+\omega \varepsilon_k}O_{\sigma(k+\omega\varepsilon_k)}^{+1}$. This follows from 
the fact that the auxiliary arcs $\zeta_k$ and $\xi_k$ intersect $\beta_k^\#$ only at the 
points $I_k^{-1}$ and $O_k^{+1}$, hence the copies of $\zeta_k$ and $\xi_k$ inside the bubble 
coming 
from $\beta$  cross it from side to side and at the same time they also cross the grafting annulus. 
As observed before, the procedure does not use  the same auxiliary object twice; this guarantees 
that coming back to the grafting annulus does not result in new intersections, so that these paths  
developing to $\widehat{\beta}'$ are actually the boundary of a new bubble. Debubbling with respect 
to this new bubble gives the desired unbranched structure $\sigma_0'$ with a bubbleable arc 
$\beta'$ such that $Bub(\sigma_0,\beta)=Bub(\sigma_0',\beta')$. Notice that by construction 
$\beta'$ 
does not intersect any real component of $\sigma_0'$, because $\widehat{\beta}'$ sits entirely in 
$\hyp^2$.
\endproof

 Depending on the intersection pattern between $\beta $ and $A_\gamma$, we have different 
possibilities for what $\sigma_0'$ looks like. We are in particular interested in the easiest case, 
which is the one in which $\beta$ crosses $A_\gamma$ just once: the structure $\sigma_0'$ of the 
previous result is exactly the one obtained by degrafting $\sigma_0$ with respect to $A_\gamma$, as 
established by the following result, which provides a converse to \cite[Theorem 5.1]{CDF}.

\begin{cor}[Degrafting Lemma]\label{degraft}
Let $\rho:\pi_1(S)\to \pslr$ be Fuchsian.
 Let $\sigma_0$ be a complex projective structure with holonomy $\rho$, $A_\gamma \subset 
\sigma_0$ a grafting annulus and $\beta \subset \sigma_0$ a bubbleable arc which 
transversely crosses $A_\gamma$ just once. Then there exist a complex projective structure 
$\sigma_0'$ with holonomy $\rho$ and a bubbleable arc $\beta'\subset \sigma_0'$ such that 
$\sigma_0=Gr(\sigma_0',\gamma)$ and $Bub(\sigma_0,\beta)=Bub(\sigma_0',\beta')$.
\end{cor}
\proof
In the previous notations, we have that $\alpha_2=\alpha_1$. Therefore the new bubble 
produced by the above procedure does a full turn around $A_\gamma$ and encompasses 
the whole real curve contained in it before leaving it. Debubbling with respect to this bubble 
produces a structure which has no real curves in the homotopy class of $\gamma$. By Goldman 
classification (see \cite[Theorem C]{GO}), it must be the structure obtained by degrafting 
$\sigma_0$ with respect to $A_\gamma$.
\endproof

\subsection{Crossing a grafting region}\label{s_region}
We now address the more general case in which $\beta$ might cross a grafting region coming from a 
multigrafting, hence we resume the notation  $A_\gamma=A^1_\gamma \cup \dots 
\cup A^M_\gamma$ for the grafting region obtained by grafting $M$ times the simple closed geodesic 
$\gamma$ of $\sigma_\rho$; recall that $A_\gamma$ is obtained by taking $M$ copies of $(\cp 
\setminus \widehat{\gamma})/\rho(\gamma)$ and gluing them in a chain along their geodesic 
boundaries, so that we see $M+1$ parallel copies of the geodesic $\gamma$.\par
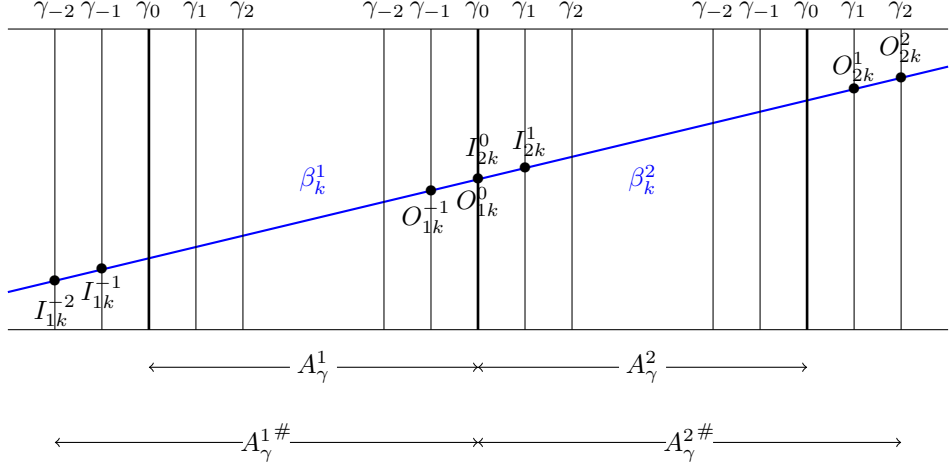
\begin{figure}[h]

\begin{center}
\begin{tikzpicture}[xscale=0.625]

 \draw (-10,-2) -- (10,-2);
 \draw (-10,2) -- (10,2);
 \draw[blue,thick] (-10,-1.5) -- (10,1.5);
 \draw[line width=1pt] (-7,2) -- (-7,-2);
 \draw[line width=1pt] (0,2) -- (0,-2);
 \draw[line width=1pt] (7,2) -- (7,-2);
\node at (-3.5,0) {\color{blue} $\beta_k^1$};
\node at (3.5,0) {\color{blue} $\beta_k^2$};

 \foreach \x in {-9,...,-5}
  \draw (\x,2) -- (\x,-2);
 \foreach \x in {-2,...,2}
  \draw (\x,2) -- (\x,-2);
 \foreach \x in {5,...,9}
  \draw (\x,2) -- (\x,-2);
  
  \node at (-9,2.25) {$\gamma_{-2}$};
  \node at (-8,2.25) {$\gamma_{-1}$};
  \node at (-7,2.25) {$\gamma_{0}$};
  \node at (-6,2.25) {$\gamma_{1}$};
  \node at (-5,2.25) {$\gamma_{2}$};
\begin{scope}[xshift=7cm]
 \node at (-9,2.25) {$\gamma_{-2}$};
  \node at (-8,2.25) {$\gamma_{-1}$};
  \node at (-7,2.25) {$\gamma_{0}$};
  \node at (-6,2.25) {$\gamma_{1}$};
  \node at (-5,2.25) {$\gamma_{2}$};
\end{scope}
\begin{scope}[xshift=14cm]
 \node at (-9,2.25) {$\gamma_{-2}$};
  \node at (-8,2.25) {$\gamma_{-1}$};
  \node at (-7,2.25) {$\gamma_{0}$};
  \node at (-6,2.25) {$\gamma_{1}$};
  \node at (-5,2.25) {$\gamma_{2}$};
\end{scope}
  
  \draw[<-] (-7,-2.5) -- (-4,-2.5);
    \node at (-3.5,-2.5) {$A_\gamma^1$};
  \draw[->] (-3,-2.5) -- (0,-2.5);
  \draw[<-] (-9,-3.5) -- (-5,-3.5);
    \node at (-4.5,-3.5) {${A_\gamma^1}^\#$};
  \draw[->] (-4,-3.5) -- (0,-3.5);
  
  \begin{scope}[xshift=7cm]
   \draw[<-] (-7,-2.5) -- (-4,-2.5);
    \node at (-3.5,-2.5) {$A_\gamma^2$};
  \draw[->] (-3,-2.5) -- (0,-2.5);
  \end{scope}
  \begin{scope}[xshift=9cm]
  \draw[<-] (-9,-3.5) -- (-5,-3.5);
    \node at (-4.5,-3.5) {${A_\gamma^2}^\#$};
  \draw[->] (-4,-3.5) -- (0,-3.5);
  \end{scope}

\node at (-9,-1.35) { $\bullet$};
\node at (-8,-1.2) { $\bullet$};
\node at (-1,-0.15) { $\bullet$};
\node at (-9,-1.75) {$I_{1k}^{-2}$};
\node at (-8,-1.5) {$I_{1k}^{-1}$};
\node at (-1.1,-0.5) {$O_{1k}^{-1}$};
\node at (0,-0.3) {$O_{1k}^{0}$};

\begin{scope}[xscale=-1,yscale=-1]
\node at (-9,-1.35) { $\bullet$};
\node at (-8,-1.2) { $\bullet$};
\node at (-1,-0.15) { $\bullet$};

\node at (-9,-1.75) {$O_{2k}^{2}$};
\node at (-8,-1.5) {$O_{2k}^{1}$};
\node at (-1.1,-0.5) {$I_{2k}^{1}$};
\node at (-0.1,-0.4) {$I_{2k}^{0}$};

\end{scope}
\node at (0,0) { $\bullet$};



\end{tikzpicture}
\end{center}
 \caption{The extended region $A_\gamma^\#$ and the auxiliary curves $\gamma_j$ in the case 
$M=2$.}\label{pic_extgraftregion}
 \end{figure}

What we want to do is to subdivide $A_\gamma$ in disjoint annular regions in such a way that we are 
able to follow the procedure described above for the case of a simple grafting inside each of them. 
The natural subdivision given by the grafting annuli $A_\gamma^h$ does not work:  
the procedure described above makes use of auxiliary curves parallel to $\gamma$ obtained by 
slightly enlarging the grafting annulus; if we did the same here we would see a lot of overlaps.
To solve this problem we consider more auxiliary curves on each side of the grafting 
geodesic, as many as the number $M$ of grafting annuli which compose the grafting region 
$A_\gamma=A^1_\gamma \cup \dots \cup A^M_\gamma$. For instance we can consider the curves
$\gamma_{\pm h} = \{ x \in \sigma_\rho \ | \ d(x,\gamma)= h \varepsilon \}$
for $h=1,\dots, M$ and an arbitrarily small $\varepsilon>0$. They clearly develop to the  
boundaries of the regions $\mathcal{N}_{h\varepsilon} (\widehat{\gamma})=\{ 
\widehat{x} \in \hyp^2 \ | \ d(\widehat{x},\widehat{\gamma})\leq h \varepsilon \}$. Recall that the 
grafting annuli $A_\gamma^h$ and $A_\gamma^{h+1}$ meet along a copy of the grafting geodesic 
$\gamma$, hence around each of these copies we have well defined copies of the curves $\gamma_j$ 
for $j=-M,\dots,M$, which we denote in the same way by a little abuse of notation; of course 
$\gamma_0$ is exactly $\gamma$ (see Figure \ref{pic_extgraftregion}).\par 
Given an oriented bubbleable properly embedded arc $\beta$ which transversely cross $A_\gamma$ from 
side to side, we can consider the crossings given by its intersections with the grafting annuli 
$A_\gamma^h$. Let us label the grafting annuli and the auxiliary curves $\gamma_{j}$ so that the 
first annulus met by $\beta$ is $A^1_\gamma$ and the first auxiliary curve is $\gamma_{-M}$. We 
obtain a doubly indexed family of crossings: $\beta_k^h$ will be the $k$-th time (with respect to 
the orientation of $\beta$) that $\beta$ crosses the annulus $A^h_\gamma$. We begin with some 
remarks. First of all the transversality assumption implies that once $\beta$ enters 
in $A_\gamma$ it has to leave on the other side, so that in each annulus $A^h_\gamma$ we see the 
same number of crossings, which we call $N$. Secondly since $\beta$ is bubbleable and all the 
grafting annuli have the same developed image, we get that the crossings 
$\beta_k^1,\dots,\beta_k^M$ have the same coherence and hence induce  the same orientation of 
$\gamma$. Therefore we can consistently orient everything using $\beta_1^1$. As before this allows 
us to order the crossings according to the cyclic order in which they appear along this 
orientation; 
once again we obtain a doubly indexed family of crossings $\alpha_k^h=\beta_{\sigma(k)}^h$ for some 
permutation $\sigma \in \mathfrak{S}_N$ such that $\sigma(1)=1$. Notice that the permutation 
$\sigma$ is the same for all the annuli $A_\gamma^1,\dots,A_\gamma^M$ because  $\beta$ is 
embedded, and that the exit point for $\beta_k^h$ coincides with the entry point of 
$\beta_k^{h+1}$.\par

Exactly as before we need to define some auxiliary points and arcs. Recall that around each 
parallel copy of $\gamma$ we have a whole package of curves which we have 
labelled $\gamma_{-M},\dots,\gamma_{M}$. 
Let us denote by $A_\gamma^\#$ the annular region containing $A_\gamma$ and bounded by $\gamma_{\pm 
M}$, and by $A_\gamma^{h\#}$ the annular region contained in $A_\gamma^\#$, bounded by 
$\gamma_{-M+2h-1\pm 1}$ and containing exactly two real curves, for $h=1,\dots,M$; roughly speaking 
these regions are obtained by slightly moving $A_\gamma^h$ by a certain amount of $\varepsilon$ 
depending on the index $h$.
Notice that the annuli $A_\gamma^{h\#}$ have disjoint interior and meet pairwise along some 
$\gamma_{j}$: more precisely $A_\gamma^{h\#}$ meets $A_\gamma^{h+1\#}$ along $\gamma_{-M+2h}$. Let 
us define the crossing $\beta_k^{h\#}=\beta \cap A_\gamma^{h\#}$ and label the intersections of 
$\beta_k^{h\#}$ with $\gamma_{-M+2h-2},\gamma_{-M+2h-1}$ and $\gamma_{-M+2h}$ by 
$I_{hk}^{-M+2h-2},I_{hk}^{-M+2h-1},O_{hk}^{-M+2h-1},O_{hk}^{-M+2h}$ in such a way that they appear 
in this order along $\beta$ (see Figure \ref{pic_extgraftregion}). Notice that 
$O_{hk}^{-M+2h}=I_{h+1,k}^{-M+2(h+1)-2}$ and that a point 
whose apex is $j$ belongs to an auxiliary curve labelled $\pm j$, according to whether that 
crossing enters the grafting region on the same side as $\beta_1$ or not.\par
Finally let us define $\widehat{\zeta}_{hk}$ to be the geodesic from $\widehat{I}_{hk}^{-M+2h-2}$ 
to $\widehat{O}_{hk}^{-M+2h-1}$ and $\widehat{\xi}_{hk}$ to be the one from 
$\widehat{I}_{hk}^{-M+2h-1}$ to $\widehat{O}_{hk}^{-M+2h}$, in complete analogy to the case of a
simple grafting. Then we apply the same procedure described in that case modifying a crossing 
$\beta_k^{h\#}$ inside the annulus $A_\gamma^{h\#}$. Notice that $A_\gamma^{h\#}$ is almost as good 
as a genuine grafting annulus, in the sense that the open annular subregion between two copies of 
$\gamma_{-M+2h-1}$ is injectively developed.

\begin{prop}\label{onegraftingregion}
Let $\rho:\pi_1(S)\to \pslr$ be  Fuchsian.
 Let $\sigma_0$ be a complex projective structure with holonomy $\rho$, $A_\gamma=A^1_\gamma \cup 
\dots \cup A^M_\gamma \subset \sigma_0$ a grafting annulus and $\beta \subset \sigma_0$ a 
bubbleable 
arc which transversely crosses $A_\gamma$ and avoids all other grafting regions of $\sigma_0$. Then 
there exist a complex projective structure $\sigma_0'$ with the same holonomy $\rho$ and a 
bubbleable arc $\beta'\subset \sigma_0'$ which avoids all the real curves of $\sigma_0'$ and such 
that $Bub(\sigma_0,\beta)=Bub(\sigma_0',\beta')$.
\end{prop}
\proof
The strategy is the same as in the case of a simple grafting (i.e. $M=1$, see 
\ref{onegraftingannulus}), with the only difference 
that the procedure which resolves the crossing $\alpha_k^h$ must take place inside the annular 
region $A_\gamma^{h\#}$. These regions are precisely defined so that what happens inside one of 
them is completely independent from what happens inside the adjacent ones.
\endproof

\subsection{Crossing several grafting regions}\label{s_general}
Now that the ideas and the main construction have been explained in detail for the case of one 
grafting annulus and one grafting region, let us handle the general case and prove the main result 
of this paper.
\begin{thm}\label{severalgraftingregions}
Let $\rho:\pi_1(S)\to \pslr$ be  Fuchsian.
 Let $\sigma_0$ be a complex projective structure with holonomy $\rho$ and $\beta \subset \sigma_0$ 
a bubbleable arc which transversely crosses any grafting region it meets. Then there exist a 
complex 
projective structure $\sigma_0'$ with the same holonomy $\rho$ and a bubbleable arc $\beta'\subset 
\sigma_0'$ which avoids all the real curves of $\sigma_0'$ and such that 
$Bub(\sigma_0,\beta)=Bub(\sigma_0',\beta')$.
\end{thm}
\proof
The strategy is to use the same technique used in \ref{onegraftingannulus} and 
\ref{onegraftingregion} in any grafting annulus or region met by $\beta$. Notice that now  between 
two crossing of a grafting region $A_\gamma$ it is possible that $\beta$  meets some other grafting 
region $A_\delta$, for a different homotopy class $\delta$. If we tried to resolve the 
intersections between $\beta$ and $A_\gamma$, it would be impossible to control the behaviour of 
the developed images of the subarcs coming from $\beta\cap A_\delta$, i.e. to prove that the above 
procedure produces an injectively developed path in $\cp$. A  way to avoid this kind of problems, 
is to  apply the procedure of \ref{onegraftingregion}  \textit{simultaneously} to all grafting 
regions met by $\beta$, without trying to handle different grafting regions one by one. 

To check that everything works as desired, it is enough to observe that for any couple of disjoint 
simple closed geodesics $\gamma$ and $\delta$ on the underlying uniformizing structure 
$\sigma_\rho$ all their developed images are disjoint; therefore not only are the associated 
grafting regions $A_\gamma$ and $A_\delta$ disjoint, but it is also possible to adjust the width of 
the extended grafting regions $A_\gamma^\#$ and $A_\delta^\#$ in such a way that the developed 
images of auxiliary arcs associated to different crossings are disjoint.
This construction realises $Bub(\sigma_0,\beta)$ as a bubbling of another structure $\sigma_0'$ 
along an arc $\beta'$ as before; by definition it avoids the real curves, exactly because we have 
replaced the portion crossing the grafting annuli with small geodesic arcs entirely contained in 
the upper-half plane.
\endproof

\subsection{The multi(de)grafting lemma}
We have already mentioned that \cite[Theorem 5.1]{CDF} states that any simple grafting can be 
obtained via a sequence of one bubbling and one debubbling, and we have proved an analogous 
statement for a simple degrafting in  \ref{degraft} under the assumption of Fuchsian 
holonomy. Under the same assumption, we can now obtain the same statement for any 
multi(de)grafting, 
by \ref{severalgraftingregions}. In particular we can show that it is possible to completely 
degraft 
a structure and recover the uniformizing structure $\sigma_\rho$ by just one bubbling and one 
debubbling.
\begin{cor}[Multi(de)grafting Lemma]\label{multidegraftviabubbltohyp}
Let $\rho:\pi_1(S)\to \pslr$ be Fuchsian and $\sigma_\rho$ the associated uniformizing 
structure.
 Let $\sigma_0$ be a complex projective structure with holonomy $\rho$ and $\beta \subset \sigma_0$ 
a bubbleable arc which transversely crosses all the grafting region of $\sigma_0$ exactly once. 
Then there exists a bubbleable arc $\beta_\rho\subset \sigma_\rho$ such that 
$Bub(\sigma_0,\beta)=Bub(\sigma_\rho,\beta_\rho)$.
\end{cor}
\proof
Let $A_{\gamma_1},\dots, A_{\gamma_n}$ be the grafting regions of $\sigma_0$. By 
\ref{severalgraftingregions} in $Bub(\sigma_0,\beta)$ we can find another bubble avoiding all real 
curves. Debubbling with respect to this bubble gives an unbranched structure 
without real curves, as in \ref{degraft}; once again by Goldman classification in 
\cite[Theorem C]{GO} it must be the uniformizing structure.
\endproof
Notice that the roles of $\sigma_0$ and $\sigma_\rho$ are symmetric in 
the above statement, in the sense that the same proof also proves that any multigrafting on 
$\sigma_\rho$ can be obtained via a sequence of just one bubbling and one debubbling. In particular 
we get the following bound on the number of moves needed to join a couple of complex projective 
structures with the same Fuchsian holonomy.
\begin{cor}\label{boundunbranched}
 Let $\rho:\pi_1(S)\to \pslr$ be Fuchsian and $\sigma,\tau$ be a couple of complex 
projective structures with holonomy $\rho$. Then it is possible to go from one to the other via 
a sequence of at most two bubblings and two debubblings.
\end{cor}
\proof
By Goldman's Theorem each of them is a multigrafting on $\sigma_\rho$. We can completely degraft 
$\sigma$ and reach $\sigma_\rho$ with one bubbling and one debubbling thanks to 
\ref{multidegraftviabubbltohyp}. Then we perform another bubbling and another debubbling to perform 
the multigrafting on $\sigma_\rho$ which produces $\tau$.
\endproof
An alternative proof can be obtained by replacing Goldman's Theorem by a result of 
Calsamiglia-Deroin-Francaviglia (see \cite[Theorem 1.1]{CDF2}) according to which any couple of 
complex projective structures with the same  Fuchsian holonomy are joined by a sequence of 
two multigraftings.

\subsection{The point of view of branched structures}
As observed in \ref{branch!}, performing a bubbling introduces a couple of simple branch points on 
the surface; therefore, properly speaking, it is not a deformation of complex projective 
structures, but of branched complex projective structures. In the previous sections these have been 
used just as a tool to study the grafting surgery on unbranched structures. Here we want to 
reformulate the main statements from the point of view of the intrinsic geometry of 
branched structures. Let us denote by $\mathcal{M}_{k,\rho}$ the moduli space of (marked) 
branched complex projective structures with a fixed holonomy $\rho$ and $k$ branch points (counted 
with multiplicity). First of all \ref{severalgraftingregions} can be restated as follows.
\begin{thm}\label{nowelldefunbranched}
 Let $\rho:\pi_1(S)\to \pslc$ be Fuchsian.
 Let $\sigma_0\in \mathcal{M}_{0,\rho}$, $\beta \subset \sigma_0$ be a bubbleable arc and $\sigma 
= Bub(\sigma_0,\beta) \in \mathcal{M}_{2,\rho}$. Assume that every time $\beta$ intersects some 
grafting region of $\sigma_0$ it actually crosses it transversely.
Then $\sigma$ is also a bubbling over some other $\sigma_0'\in \mathcal{M}_{0,\rho}$ 
along a bubbleable arc $\beta' \subset \sigma_0'$ which avoids the real curves of $\sigma_0'$.
\end{thm}
In a previous paper of the author (see \cite[Theorem 5.9]{R}) it was proved that there exists a 
connected open dense subspace of $\mathcal{M}_{2,\rho}$ consisting of structures obtained by 
bubbling unbranched structures. A consequence of the previous statement is that these structures do 
not have such a thing as a canonical underlying unbranched structure in general, in the sense that 
the same branched structure may arise as a bubbling over different unbranched structures. We 
observed there (see \cite[§5.1]{R}) that it is generically possible to join any couple of structures 
in $\mathcal{M}_{0,\rho}\cup \mathcal{M}_{2,\rho}$ by a finite sequence of bubblings and 
debubblings. The results of this paper allow us to give a uniform upper bound on the length of such 
a sequence.

\begin{cor}\label{6steps}
 Let $\rho:\pi_1(S)\to \pslc$ be Fuchsian. There is a connected, open and dense subspace 
$\mathcal{B}\subset \mathcal{M}_{2,\rho}$ such that if  $\sigma,\tau \in \mathcal{M}_{0,\rho}\cup 
\mathcal{B}$ then $\tau$ is obtained from $\sigma$ via a sequence 
of at most three bubblings and three debubblings.
\end{cor}
\proof
If $\sigma$ and $\tau$ are unbranched then this follows directly from \ref{boundunbranched} 
above. Otherwise we can take  $\mathcal{B}$ to be the space of structures obtained by bubbling 
unbranched structures provided by \cite[Theorem 5.9]{R}. If $\sigma,\tau \in \mathcal B$, then by 
definition of $\mathcal B$ we can find two unbranched structures $\sigma_0,\tau_0 \in \mathcal 
M_{0,\rho}$ such that $\sigma$ and $\tau$ are obtained by bubbling $\sigma_0$ and $\tau_0$ 
respectively, and then apply \ref{boundunbranched} to them.
\endproof

\begin{figure}[h]

\begin{center}
\begin{tikzpicture}
\node at (-6,0) {$\sigma$};

\draw[->] (-5.5,0) to[out=0,in=180] (-3.5,0);
\node at (-4.5,.5) {$1$ debub};

\node at (-3,0) {$\sigma_0$};

\draw[->] (-2.5,0) to[out=0,in=180] (-.5,0);
\node at (-1.5,.5) {$1$ bub};
\node at (-1.5,-.5) {$1$ debub};

\node at (0,0) {$\sigma_\rho$};

\draw[->] (.5,0) to[out=0,in=180] (2.5,0);
\node at (1.5,.5) {$1$ bub};
\node at (1.5,-.5) {$1$ debub};

\node at (3,0) {$\tau_0$};

\draw[->] (3.5,0) to[out=0,in=180] (5.5,0);
\node at (4.5,.5) {$1$ bub};

\node at (6,0) {$\tau$};

\end{tikzpicture}
\end{center}
\end{figure}

\printbibliography

\end{document}